\RequirePackage{fix-cm}
\documentclass[smallextended]{svjour3}       
\smartqed  
\usepackage{
amsmath,
amssymb,
hyperref,
mathrsfs,
graphicx,
stmaryrd,
amscd,
 amsfonts, stmaryrd,latexsym, xargs}
\usepackage{color}

\usepackage[colorinlistoftodos]{todonotes}
 \presetkeys{todonotes}%
{inline,backgroundcolor=gray!20,bordercolor=gray!30}{}
\tikzset{/tikz/notestyleraw/.append style={text=black}}

\newcommand\todoin[2][]{\todo[inline, caption={2do}, #1]{
 \begin{minipage}{\textwidth-3pt}#2\end{minipage}
 }}

\newtheorem{thm}{Theorem}[section]

\newtheorem{lem}[thm]{Lemma}
\newtheorem{defn}[thm]{Definition}

\newtheorem{prop}[thm]{Proposition}

\newtheorem{rmk}[thm]{Remark}
\newcommand{\beq}{\begin{equation}}
\newcommand{\eeq}{\end{equation}}
\newcommand{\be}{\begin{eqnarray}}
\newcommand{\ee}{\end{eqnarray}}
\newcommand{\ben}{\begin{eqnarray*}}
\newcommand{\een}{\end{eqnarray*}}

\newcommand{\beal}{\begin{aligned}}
\newcommand{\enal}{\end{aligned}}
\newcommand{\bc}{\begin{cases}}
\newcommand{\ec}{\end{cases}}

\newcommand{\R}{\mathbb{R}}
\newcommand{\PP}{\mathbb{P}}

\newcommand{\N}{\mathbb{N}}

\newcommand{\cG}{\mathcal{G}}

\newcommand{\cA}{\mathcal{A}}
\newcommand{\cO}{\mathcal{O}}
\newcommand{\cM}{\mathcal{M}}

\newcommand{\wt}{\widetilde}
\newcommand{\wh}{\widehat}
\newcommand{\ol}{\overline}

\newcommand{\om}{\omega}
\newcommand{\lb}{\lambda}
\newcommand{\Lb}{\Lambda}

\begin{document}

\title{On the negative limit of viscosity solutions for discounted Hamilton-Jacobi equations}
\subtitle{}

\titlerunning{Negative limit of viscosity solutions for discounted HJ equations}        

\author{Ya-Nan Wang
\and Jun Yan
\and
Jianlu Zhang
}


\institute{Ya-Nan Wang \at
School of Mathematical Sciences, Nanjing Normal University, Nanjing, 210097, China\\
\email{yananwang@fudan.edu.cn}\\
Jun Yan \at
School of Mathematical Sciences, Fudan University and Shanghai Key Laboratory
for Contemporary Applied Mathematics, Shanghai 200433, China\\
\email{yanjun@fudan.edu.cn}\\
Jianlu ZHANG \at
              Hua Loo-Keng Key Laboratory of Mathematics \& Mathematics Institute, Academy of Mathematics and systems science, Chinese Academy of Sciences, Beijing 100190, China \\
              \email{jellychung1987@gmail.com}           
}

\date{Received: date / Accepted: date}

\maketitle

\begin{abstract}
Suppose  $M$ is a closed Riemannian manifold. For a $C^2$ generic  (in the sense of Ma\~n\'e) Tonelli Hamiltonian $H: T^*M\rightarrow\R$, the minimal viscosity solution $u_\lb^-:M\rightarrow \R$
of
the negative discounted equation
\[
-\lb u+H(x,d_xu)=c(H),\quad x\in M,\ \lb>0
\]
with the Ma\~n\'e's critical value $c(H)$ 
 converges to a uniquely established viscosity solution $u_0^-$ of the critical Hamilton-Jacobi equation
\[
H(x,d_x u)=c(H),\quad x\in M
\]
as $\lb\rightarrow 0_+$. We also propose a dynamical interpretation of $u_0^-$.


\end{abstract}

\tableofcontents

\section{Introduction}\label{s1}
Suppose $M$ is a closed manifold equipped with a smooth Riemannian metric $g$. $TM$ (resp. $T^*M$) denotes the tangent (resp. cotangent) bundles on $M$, on which a coordinate $(x,v) \in TM$ (resp. $(x,p) \in T^*M$) is involved. Since  $M$ is compact, we can use $| \cdot |_x$ as the norm induced by $g$ on both the fiber $T_xM$ and  $T_x^*M$, with a slight abuse of notations. A $C^2$ function $H: T^*M\times\R\to\R$ is called a {\sf Tonelli Hamiltonian}, if
\begin{enumerate}
    \item {\sf (Positive definiteness)} $\forall\ (x,p)\in T^*M$, the Hessian matrix
    $\partial_{pp}H(x,p)$ is positive definite;
    \item  {\sf (Superlinearity)} $H(x, p)$ is superlinear in the fibers, i.e.,
	\[ \lim_{|p|_x\to \infty} H(x, p)/|p|_x=+\infty,\quad \forall  x\in M.\] 
\end{enumerate}
For such a Hamiltonian, we consider the {\sf negative discounted equation} 
\beq\label{eq:dis}\tag{HJ$_\lb^-$}
-\lb u+H(x,d_xu)=c(H),\quad x\in M,\ \lb>0
\eeq
with the {\sf Ma\~n\'e critical value} (see \cite{CIPP})
\[
c(H):=\inf_{u\in C^\infty(M,\R)}\sup_{x\in M} H(x,d_xu).
\]
We can prove that the viscosity solution of (\ref{eq:dis}) indeed exists but unnecessarily unique (see Sec. \ref{s2}).
Nonetheless, 
the {\sf ground state solution} defined by
\[
u_\lb^-:=\inf\{u(x)|x\in M, u \mbox{ is a viscosity solution of (\ref{eq:dis}) }\}
\]
is uniquely identified and satisfies a bunch of fine properties (see Sec. \ref{s3} for details). That urges us to explore the convergence of $u_\lb^-$ as $\lb\rightarrow 0_+$:
\begin{thm}[Main 1]\label{thm:1}
For a $C^2$ smooth Tonelli Hamiltonian $H(x,p)$ generic in the sensee of Ma\~n\'e\footnote{See Definition \ref{defn:man-gene} for the definitions of Ma\~n\'e's genericity. 
}, 
the ground state solution $u_\lb^-$ of \eqref{eq:dis} converges to a specified viscosity solution $u_0^-$ of 
\beq\label{eq:hj}\tag{$\text{HJ}_0$}
H(x,d_x u)=c(H),\quad x\in M
\eeq
 as $\lb\rightarrow 0_+$. 
\end{thm}
\begin{rmk}
Actually, we proved the convergence of $u_\lb^-$ for any $C^2$ Tonelli Hamiltonian $H(x,p)$ satisfying the following:
\todoin{
for any sequence $\lb_n\rightarrow 0_+$ as $n\rightarrow+\infty$,
\beq\label{upper-aubry}\tag{$\star$}
\Big(\bigcap_{N>0}\ol{\bigcup_{n\geq N}\cG_{\lb_n}}\Big)\bigcap\Lb_i\neq\emptyset, \quad \forall \Lb_i\in \cA\slash d_c.
\eeq
}
\noindent see Sec. \ref{s2} for the definition of $\cG_\lb$ and Appendix \ref{a1} for $\cA$, $\cA\slash d_c$. This is a generic condition for $C^2$ Tonelli Hamiltonian $H(x, p)$.\medskip


By using a variational analysis of contact Hamiltonian systems developped in \cite{WWYNON,WWYJMPA,WWY,WWY2}, a breakthrough on the vanishing discount problem in the negative direction was recently made in \cite{DW}. Precisely,  for any $C^3$ smooth Tonelli Hamiltonian $H(x,p)$ satisfying 
\todoin{
 constant functions are subsolutions of \eqref{eq:hj},\hspace{112pt}$(\lozenge)$
 }
\noindent  they proved the convergence of the ground state solution $u_\lb^-$ of (\ref{eq:dis}) as $\lb\rightarrow 0_+$.
 Actually, any $H(x,p)$ satisfying $(\lozenge)$ has to satisfy \eqref{upper-aubry}, since for any $\lb>0$ the associated $\cA\subset \cG_\lb$ (proved in Proposition 3.2 of \cite{DW}). However, $(\lozenge)$ is not generic, that's the reason we find a more general substitute  $(\star)$ in this paper. \medskip

 The significance of the convergence of viscosity solutions for discounted equations was firstly proposed by Lions, Papanicolaou and Varadhan in 1987 \cite{LPV}. In \cite{DFIZ}, a rigorous proof for the convergence was finally given for the positive discounted systems. Also in recent works \cite{CCIZ,WYZ,Z}, the convergence of solutions for generalized $1^{st}$ order PDE was discussed. Comparing to these works, the negative discount in \eqref{eq:dis} brings new difficulties to prove the existence of viscosity solutions, not to mention the convergence. By using a dual Lagrangian approach, we reveal that the negative discount limit of solutions actually conjugates to the positive discounted limit, see Sec. \ref{s3} for details.
\end{rmk}

\vspace{10pt}

\noindent{\bf Organization of the article.} In Sec. \ref{s2}, we discuss the convergence of the Lax-Oleinik semigroup for \eqref{eq:dis} and the properties of the ground state solution $u_\lb^-$. In Sec. \ref{s3}  we finish the proof of Theorem \ref{thm:1}.  For the consistency and readability, a brief review of the Aubry-Mather theory and some properties of the Lax-Oleinik semigroup are moved to the Appendix.

\vspace{5pt}

\noindent{\bf Acknowledgement.} The first author is supported by National Natural Science Foundation of China (Grant No.11501437). The second author is supported by National Natural Science Foundation of China (Grant No. 11631006 and 11790272) and Shanghai Science and Technology Commission (Grant No. 17XD1400500). The third author is supported by the National Natural Science Foundation of China (Grant No. 11901560). 

\section{Generalities}\label{s2}



%

\begin{defn}[weak KAM solution]
For $\lb\geq0$, a function $u:M\rightarrow \R$ is called a {\sf backward (resp. forward) $\lb-$weak KAM solution} if it satisfies:
\begin{itemize}
\item $u\prec_\lb L+c(H)$, i.e. for any $(x,y)\in M\times M$ and $a\leq b$, we have 
\[
e^{-\lb b}u(y)-e^{-\lb a}u(x)\leq h_\lb^{a,b}(x,y)
\]
with 
\be\label{eq:sta}
 h_\lb^{a,b}(x,y):=\inf_{\substack{\gamma\in C^{ac}([a,b],M)\\\gamma(a)=x,\gamma(b)=y}}\int_a^be^{-\lb s}\Big(L(\gamma,\dot\gamma)+c(H)\Big)ds
\ee
\item for any $x\in M$ there exists a curve $\gamma_{x,\lb}^-:(-\infty,0]\rightarrow M$ (resp. $\gamma_{x,\lb}^+:[0,+\infty)\rightarrow M$) ending with (resp. starting from) $x$, such that for any $s<t\leq 0$ (resp. $0\leq s<t$), 
\[
e^{-\lb t}u(\gamma_{x,\lb}^-(t))-e^{-\lb s}u(\gamma_{x,\lb}^-(s))=\int_s^te^{-\lb\tau}\big(L(\gamma_{x,\lb}^-,\dot\gamma_{x,\lb}^-)+c(H)\big)d\tau.
\]
\[
\bigg(resp.\quad e^{-\lb t}u(\gamma_{x,\lb}^+(t))-e^{-\lb s}u(\gamma_{x,\lb}^+(s))=\int_s^te^{-\lb \tau}\big(L(\gamma_{x,\lb}^+,\dot\gamma_{x,\lb}^+)+c(H)\big)d\tau.\bigg)
\]
Such a $\gamma_{x,\lb}^-$ (resp. $\gamma_{x,\lb}^+$) is called a {\sf backward (resp. forward) calibrated curve} by $u$.
\end{itemize} 
\end{defn}


Notice that
\beq\label{eq:for-weak}\tag{S$^+$}
{u}_\lb^+(x):=-\inf_{\substack{\gamma\in C^{ac}([0,+\infty),M)\\\gamma(0)=x}}\int^{+\infty}_0e^{-\lb \tau}\Big(L(\gamma(\tau),\dot\gamma(\tau))+c(H)\Big) d\tau
\eeq
is the unique forward $\lb-$weak KAM solution of \eqref{eq:dis} (see Proposition \ref{prop:dom-fun} for details). However, it's not straightforward to get one backward $\lb-$weak KAM solution for \eqref{eq:dis}. Usually, that relies on the {\sf  backward Lax-Oleinik operator}  
\[
T_t^{\lb,-}: C(M,\R)\rightarrow C(M,\R)
\]
via
\[
T_t^{\lb,-}\phi(x)=\inf_{\gamma(0)=x}\bigg\{e^{\lb t}\phi(\gamma(-t))+\int^0_{-t} e^{-\lb \tau}\Big(L(\gamma(\tau),\dot\gamma(\tau))+c(H)\Big)d\tau\bigg\},\quad t\geq 0.
\]
If we can find a fixed point of $T_t^{\lb,-}$, then it automatically becomes a backward $\lb-$weak KAM solution of \eqref{eq:dis}, then becomes a viscosity solution of \eqref{eq:dis} (see (6) \& (7) of Lemma \ref{lem2:proT}).\\

Different from the approach in \cite{DW,WWY,WWY2}, here we proved the same conclusions as \cite{DW} in our own way:

\begin{prop}[ground state]\label{prop:dom-fun}
For $\lb>0$, 
\be
u_\lb^-:=\lim_{t\rightarrow+\infty}T_t^{\lb,-}u_\lb^+
\ee
 exists as a backward $\lb-$weak KAM solution of \eqref{eq:dis} and satisfies the following:
\begin{enumerate}
\item [(a)] $u_\lb^-(x)$ is a viscosity solution of \eqref{eq:dis};
\item [(b)] $u_\lb^-(x)\geq u_\lb^+(x)$ for all $x\in M$ and the equality holds on and only on the following set \footnote{the following curve $\gamma$ is called {\sf  globally calibrated} by $u_\lb^+$}
\be\label{defn:g-cali}
\cG_\lb&:=&\bigg\{x\in M\Big|\exists\  \gamma:\R \rightarrow M \text{with } \gamma(0)=x, \  e^{-\lambda b}u_\lb^+(\gamma(b))-e^{-\lambda a}u_\lb^+(\gamma(a))\nonumber\\
&=&\int^b_ae^{-\lambda\tau}(L(\gamma(\tau),\dot{\gamma}(\tau))+c(H))\mbox{d}\tau,\ \forall \ a<b\in\R\bigg\}.
\ee
\item [(c)] $u_\lb^-$ is the minimal viscosity solution of \eqref{eq:dis}. Namely, any fixed point of $T_t^{\lb,-}:C(M,\R)\rightarrow C(M,\R)$ has to be greater than $u_\lb^-$.
\end{enumerate}
\end{prop}
\proof

As we can see,
\be\label{eq:tran}
\wh u_\lb^-(x)&:=&-u_\lb^+(x)\nonumber\\
&=&\inf_{\substack{\gamma\in C^{ac}([0,+\infty),M)\\\gamma(0)=x}}\int^{+\infty}_0e^{-\lb t}\Big(L(\gamma(t),\dot\gamma(t))+c(H)\Big) dt\nonumber\\
&=&\inf_{\substack{\gamma\in C^{ac}((-\infty,0],M)\\\gamma(0)=x}}\int_{-\infty}^0e^{\lb t}\Big(
\wh L(\gamma(t),\dot\gamma(t))+c(H)\Big) dt,
\ee
where $\wh L(x,v):=L(x,-v)$ for $(x,v)\in TM$. Due to Appendix 2 in \cite{DFIZ}, $\wh u_\lb^-$ is the unique viscosity solution of 
\be\label{eq:sym-dis}
\lb u+\wh H(x,d_x u)=c(H),\quad x\in M,\lb>0
\ee
with $\wh H(x,p):=H(x,-p)$ for $(x,p)\in T^*M$. Besides, for any $x\in M$ the curve $\wh\gamma_{x,\lb}^-:t\in(-\infty,0]\rightarrow M$ ending with $x$ and attaining the infimum of \eqref{eq:tran} always exists, and $\|\dot{\wh\gamma}_{x,\lb}^-(t)\|_{L^\infty}\leq\alpha$ for some constant $\alpha>0$ independent of $x$ and $\lb$. Therefore, $\gamma_{x,\lb}^+(t):=\wh\gamma_{x,\lb}^-(-t)$ defines a curve for $t\in [0,+\infty)$ starting from $x$, which is forward calibrated by $u_\lb^+$. On the other side, we can easily verify that $u_\lb^+\prec_\lb L+c(H)$, so $u_\lb^+$ in \eqref{eq:for-weak} is indeed a forward $\lb-$weak KAM solution of \eqref{eq:dis}.

Since now $\wh L(x,v)$ is $C^2-$smooth, due to the {Tonelli Theorem} and {\sf Weierstrass Theorem} in \cite{Mat}, $\gamma^+_{x,\lambda}(\tau)$ is $C^2-$smooth and $(\gamma^+_{x,\lambda}(\tau),\dot \gamma^+_{x,\lambda}(\tau))$solves the Euler-Lagrangian equation:
\begin{equation}\label{eq2:E-L}
\left\{
\begin{aligned}
&\dot x=v\\
	&\frac{\mbox{d}}{\mbox{d}t}\bigg(\frac{\partial L}{\partial v}(x,v)\bigg)-\lambda\frac{\partial L}{\partial v}(x,v)=\frac{\partial L}{\partial x}(x,v),
	\end{aligned}
	\right.
\end{equation}
for $\tau\in(0,+\infty)$. 
Hence,
$|\dot{\gamma}^+_{x,\lambda}(\tau)|\leq \alpha$ for all $\tau\in(0,+\infty)$. 
{Moreover, due to Proposition 2.5 in \cite{DFIZ}, $\{\wh u^-_\lambda(x)\}_{\lb>0}$ are equi-Lipschitz and equi-bounded, then consequently so are $\{u^+_\lambda(x)\}_{\lb>0}$.
}

Let $\Phi^s_{L,\lambda}:TM\rightarrow TM$ be the flow of \eqref{eq2:E-L}. By an analogy of \cite{MS}, we define 
\begin{align*}
	\widetilde{\Sigma}_\lb=\{(x,v)\in TM|\gamma(s)&=\pi\circ\Phi^s_{L,\lambda}(x,v)\\
	&\mbox{ is a forward calibrated curve by }u^+_{\lambda}\}
\end{align*}
where $\pi:TM\rightarrow M$ is the standard projection. Apparently, 
\begin{description}
	\item [(i)] $\widetilde{\Sigma}_\lambda$ is compact.
	\item [(ii)] $\widetilde{\Sigma}_\lambda$ is forward invariant under $\Phi^t_{L,\lambda}$, i.e.,
	$\Phi^t_{L,\lambda}(\widetilde{\Sigma}_\lambda)\subset\widetilde{\Sigma}_\lambda$, for all $t>0$. 
\end{description}	
Then
\begin{equation}
	\widetilde{\mathcal{G}}_\lb=\bigcap_{t\geq 0}\Phi^t_{L,\lambda}\big(\widetilde{\Sigma}_\lambda\big)
\end{equation}
is nonempty, compact and $\Phi_{L,\lb}^t-$invariant. Due to (a.3) of Sec. 4 in \cite{MS}, $\cG_\lb:=\pi{\wt\cG}_\lb$ has an equivalent definition as in \eqref{defn:g-cali}.\medskip


%

Next, we claim two conclusions:
\begin{itemize}
\item [(iii)] $T^{\lambda,-}_tu^+_{\lambda}(x)\geq u^+_{\lambda}(x)$, for all $x\in M, \ t>0$ (proved in Lemma \ref{lem2:Tmono});

\item [(iv)]  For any $\lambda>0$ fixed, the set $T^{\lambda,-}_t u^+_{\lambda}(x)$ are uniformly bounded from above for all $t\geq1$ and $x\in M$ (proved in Lemma \ref{lem2:Tbounded}).
\end{itemize}
Benefiting from these conclusions and (5) of Lemmas \ref{lem2:proT},
$T^{\lambda,-}_tu^+_\lambda(x)$ converges uniformly to a function (denoted by $u^-_\lambda(x))$ on $M$ as $t$ tending to $+\infty$. Due to (1) of Lemma \ref{lem2:proT}, for each $t_1>0$ we derive
\begin{align*}
	u^-_\lambda(x)&=\lim_{t\to+\infty}T^{\lambda,-}_{t_1+t}u^+_\lambda(x)
=\lim_{t\to+\infty}T^{\lambda,-}_{t_1}\circ T^{\lambda,-}_{t}u^+_\lambda(x)\\
&=
T^{\lambda,-}_{t_1}(\lim_{t\to+\infty}T^{\lambda,-}_tu^+_\lambda(x))=T^{\lambda,-}_{t_1}u^-_\lambda(x),
\end{align*}
then by (6) of Lemma \ref{lem2:proT}, $u^-_\lambda(x)$ is a backward $\lambda$-weak KAM solution of \eqref{eq:dis}, and by (7) of Lemma \ref{lem2:proT}, $u^-_\lambda(x)$ is a viscosity solution of \eqref{eq:dis}. So (a) of Proposition \ref{prop:dom-fun} is proved.\medskip

Due to (iii), we get $u_\lb^-\geq u_\lb^+$. For any $x\in\cG_\lb$, there exists a globally calibrated curve $\gamma_x:\R\rightarrow M$ passing it and achieving \eqref{eq:sta} on any interval $[a,b]\subset\R$. So  for any $t>0$, 
\begin{align*}
	T^{\lambda,-}_tu^+_\lambda(x)&\leq e^{\lambda t}u^+_\lambda(\gamma_{x}(-t))+\int^0_{-t}e^{-\lambda \tau}(L(\gamma_{x}(\tau),\dot{\gamma}_{x}(\tau))+c(H))\mbox{d}\tau\\
	&=u^+_\lambda(x),
\end{align*}
which implies $u^-_\lambda(x)= u^+_\lambda(x)$ on $\mathcal{G}_\lambda$.

On the other side, if $u^+_\lambda(z)=u^-_\lambda(z)$ for some point $z\in M$, there exists a backward calibrated curve $\gamma^-_{z,\lb}:(-\infty,0]\to M$ ending with $z$ and a forward calibrated curve $\gamma^+_{z,\lb}:[0,+\infty)\to M$ starting with $z$. Note that $u^+_\lambda(x)\leq u^-_\lambda(x)$ and $u_\lb^\pm\prec_\lambda L+c(H)$, then for any $t>0$,
\begin{align*}
	e^{-\lambda t}u^+_\lambda(\gamma^+_{z,\lambda}(t))-u^+_\lambda(z)&\leq e^{-\lambda t}u^-_\lambda(\gamma^+_{z,\lambda}(t))-u^+_\lambda(z)\\
&\leq\int^t_0e^{-\lambda \tau}(L(\gamma^+_{z,\lambda}(\tau),\dot{\gamma}^+_{z,\lambda}(\tau))+c(H))\mbox{d}\tau\\
&= e^{-\lambda t}u^+_\lambda(\gamma^+_{z,\lambda}(t))-u^+_\lambda(z),
\end{align*}
which implies $u^+_\lambda(\gamma^+_{z,\lambda}(t))=u^-_\lambda(\gamma^+_{z,\lambda}(t))$ for all $t>0$.
Similarly  $u^+_\lambda(\gamma^-_{z,\lambda}(-t))=u^-_\lambda(\gamma^-_{z,\lambda}(-t)),t>0$.
Hence, for any $a<b\in\mathbb{R}$
$$
e^{-\lambda b}u^\pm_\lambda(\gamma_{z,\lambda}(b))-e^{-\lambda a}u^\pm_\lambda(\gamma_{z,\lambda}(a))=\int^{b}_ae^{-\lambda\tau}(L(\gamma_{z,\lambda}(\tau),\dot{\gamma}_{z,\lambda}(\tau))+c(H))\mbox{d}\tau,
$$
where $\gamma_{z,\lambda}:(-\infty,+\infty)\to M$ is defined by
$$
\gamma_{z,\lambda}(\tau)
=\begin{cases}
	\gamma_{z,\lambda}^-(\tau),\tau\in (-\infty,0],\\
	\gamma_{z,\lambda}^+(\tau),\tau\in [0,+\infty).
\end{cases}
$$
So $\gamma_{z,\lambda}:\R\rightarrow M$ is globally calibrated by $u_\lb^-$ (also by $u_\lb^+$) and then $z\in\mathcal{G}_\lambda$. Finally $u_\lb^-(x)=u_\lb^+(x)$ iff $x\in\cG_\lb$, which implies item (b) of Proposition \ref{prop:dom-fun}.\medskip

To prove (c) of Proposition \ref{prop:dom-fun}, it suffices to show any backward $\lambda$-weak KAM solution $v^-$ of \eqref{eq:dis} is greater than 
$u^+_\lambda(x)$ (due to (6) and (7) of Lemma \ref{lem2:proT}). If so, 
\[
v^-(x)=T^{\lambda,-}_tv^-(x)\geq T^{\lambda,-}_tu^+_\lambda(x)
\]
for all $t>0$, then $v^-(x)\geq \lim_{t\to+\infty}T^{\lambda,-}_tu^+_\lambda(x)=u^-_\lambda(x)$.

 For each $\gamma:[0,t]\to M$ with $\gamma(0)=x$ and $t>0$, we define $\widetilde{\gamma}:[-(t+1),0]\to M$ by 
$$
\wt{\gamma}(s)=
\begin{cases}
	\gamma(s+t+1), s\in[-(t+1),-1],\\
	\beta(s), s\in [-1,0],
\end{cases}
$$
where $\beta$ is a geodesic satisfying $\beta(-1)=\gamma(t),\beta(0)=x$, and $|\dot{\beta}(s)|\leq \mbox{diam(M)}=:k_0$.
Then,
\begin{align*}
	v^-(x)&=T^{\lambda,-}_{t+1}v^-(x)\leq e^{\lambda (t+1)}v^-(x)+\int^{0}_{-(t+1)}e^{-\lambda\tau}(L(\tilde{\gamma}(\tau),\dot{\tilde{\gamma}}(\tau))+c(H))\mbox{d}\tau\\
	&\leq e^{\lambda(t+1)}v^-(x)+e^{\lambda(t+1)}\int^{t}_0e^{-\lambda\tau}(L(\gamma(\tau),\dot{\gamma}(\tau))+c(H))\mbox{d}\tau\\
	&\hspace{2.4cm}+\int^0_{-1}e^{-\lambda\tau}(L(\beta(\tau),\dot{\beta}(\tau))+c(H))\mbox{d}\tau.
\end{align*}
Hence,
\begin{align*}
	v^-(x)&\geq e^{-\lambda(t+1)}v^-(x)-\int^{t}_0e^{-\lambda\tau}(L(\gamma(\tau),\dot{\gamma}(\tau))+c(H))\mbox{d}\tau\\
	&-e^{-\lambda(t+1)}\int^0_{-1}e^{-\lambda\tau}(L(\beta(\tau),\dot{\beta}(\tau))+c(H))\mbox{d}\tau.
\end{align*}
We derive 
\[
v^-(x)\geq e^{-\lambda(t+1)}v^-(x)+A_t(x)-e^{-\lambda(t+1)}\frac{(C_{k_0}+c(H))(e^\lambda-1)}{\lambda},
\]
where $C_{k_0}=\max\{L(x,v)|(x,v)\in TM,|v|\leq k_0\}$ and
$$
A_t(x)=-\inf_{\substack{\gamma\in C^{ac}([0,t],M)\\\gamma(0)=x}}\int^t_0e^{-\lambda\tau}(L(\gamma(\tau),\dot{\gamma}(\tau))+c(H))\mbox{d}\tau.
$$
To show $v^-(x)\geq u^+_\lambda(x)$, it suffices to show $A_t(x)$ converges uniformly to $u^+_\lambda(x)$ on $M$, as $t\to+\infty$. For this purpose, we take $\gamma:[0,t]\to M$ be a minimizer of $A_t(x)$ and define $\wt{\gamma}:[0,+\infty)\to M$ by
$$
\wt{\gamma}(s)=
\begin{cases}
	\gamma(s),s\in[0,t],\\
	\gamma(t),s\in(t,+\infty).
\end{cases}
$$
Then, there exists a constant $C_0$ greater than $\max_{x\in M}|L(x,0)|$ such that 
\begin{align*}
	A_t(x)-u^+_\lambda(x)&\leq \int^{+\infty}_0e^{-\lambda\tau}(L(\tilde{\gamma}(\tau),\dot{\tilde{\gamma}}(\tau))+c(H))\mbox{d}\tau\\
	&-\int^t_0e^{-\lambda\tau}(L(\gamma(\tau),\dot{\gamma}(\tau))+c(H))\mbox{d}\tau\\
	&=\int^{+\infty}_te^{-\lambda\tau}(L(\gamma(\tau),0)+c(H))\mbox{d}\tau\\
	&\leq(C_0+c(H)) \int^{+\infty}_te^{-\lambda\tau}\mbox{d}\tau=\frac{(C_0+c(H))e^{-\lambda t}}{\lambda}.
\end{align*}
On the other hand, let $\gamma_{x,\lambda}^+:[0,+\infty)\to M$ be a minimizer of $u^+_\lambda(x)$. Recall that $\gamma_{x,\lb}^+(t)=\wh\gamma_{x,\lb}^-(-t)$, by a similar way, 
\begin{align*}
	u^+_\lambda(x)-A_t(x)&\leq -\int^{+\infty}_te^{-\lambda\tau}(L(\gamma_{x,\lambda}^+(\tau),\dot{\gamma}_{x,\lambda}^+(\tau))+c(H))\mbox{d}\tau\\
	&\leq (C(0)-c(H))\int^{+\infty}_te^{-\lambda\tau}\mbox{d}\tau=\frac{(C(0)-c(H))e^{-\lambda t}}{\lambda},
\end{align*}
where $C(0)=-\min\{L(x,v)|(x,v)\in TM\}$.
Hence, we derive $A_t(x)$ converges uniformly to $u^+_\lambda(x)$ on $M$, as $t$ tending to $+\infty$.
\qed

\begin{lem}\label{lem2:Tmono}
	$$
	T^{\lambda,-}_tu^+_{\lambda}(x)\geq u^+_{\lambda}(x),\forall x\in M, t>0.
	$$
	Hence, $T^{\lambda,-}_{t_2}u^+_\lambda(x)\geq T^{\lambda,-}_{t_1}u^+_\lambda(x)$, for $t_2\geq t_1>0$.
\end{lem}
\proof
Let  $\gamma\in C^{ac}([-t,0],M)$ with $\gamma(0)=x$ and $\eta\in C^{ac}([0,+\infty), M)$ with $\eta(0)=x$. Define $\tilde{\gamma}:[-t,+\infty)\to M$ by
$$
\tilde{\gamma}(\tau)=
\begin{cases}
	\gamma(\tau),\tau\in[-t,0],\\
	\eta(\tau),\tau\in[0,+\infty).
\end{cases}
$$
Then,
\begin{align*}
	&e^{\lambda t}u^+_{\lambda}(\gamma(-t))+\int^0_{-t}e^{-\lambda\tau}(L(\gamma(\tau),\dot{\gamma}(\tau))+c(H))\mbox{d}\tau\\
	&\geq -\int^{+\infty}_{-t}e^{-\lambda \tau}(L(\tilde{\gamma}(\tau),\dot{\tilde{\gamma}}(\tau))+c(H))\mbox{d}\tau+\int^0_{-t}e^{-\lambda\tau}(L(\gamma(\tau),\dot{\gamma}(\tau))+c(H))\mbox{d}\tau\\
	&=-\int^{+\infty}_0e^{-\lambda\tau}(L(\eta(\tau),\dot{\eta}(\tau))+c(H))\mbox{d}\tau.
\end{align*}
By the arbitrariness of $\gamma$ and $\eta$, we derive that
$
T^{\lambda,-}_t u^+_{\lambda}(x)\geq u^+_{\lambda}(x).
$
\qed
\begin{lem}\label{lem2:Tbounded}
  Given $\lambda>0$, the set $\{T^{\lambda,-}_t u^+_{\lambda}(x)|t\geq1,x\in M\}$ is bounded from above.
\end{lem}
\proof
Let $z\in \mathcal{G}_\lambda,\ \gamma_{z,\lambda}(s)$ be the globally calibrated curve by $u^+_{\lambda}$. Let $\beta:[-1,0]\to M$ be a geodesic satisfying $\beta(-1)=\gamma_{z,\lambda}(-1),\beta(0)=x$, and $|\dot{\beta}(\tau)|\leq \mbox{diam}(M)=k_0$.
Then,
\begin{align*}
	T^{\lambda,-}_tu^+_{\lambda}(x)&\leq u^+_{\lambda}(\gamma_{z,\lambda}(-t))e^{\lambda t}+\int^{-1}_{-t}e^{-\lambda\tau}(L(\gamma_{z,\lambda}(\tau),\dot{\gamma}_{z,\lambda}(\tau))+c(H))\mbox{d}\tau\\
	&+\int^0_{-1}e^{-\lambda\tau}(L(\beta(\tau),\dot{\beta}(\tau))+c(H))\mbox{d}\tau\\
	&=e^{\lambda}u^+_\lambda(\gamma_{x,\lambda}(-1))+\int^0_{-1}e^{-\lambda\tau}(L(\beta(\tau),\dot{\beta}(\tau))+c(H))\mbox{d}\tau\\
	&\leq e^{\lambda}u^+_\lambda(\gamma_{x,\lambda}(-1))+(C_{k_0}+c(H))\int^0_{-1}e^{-\lambda\tau}\mbox{d}\tau\\
	&\leq (K+\frac{C_{k_0}+c(H)}{\lambda})e^{\lambda},
\end{align*}
where $K=\|u^+_\lambda\|$.
\qed

\begin{prop}\label{prop:uni}
 $\{u_\lb^-\}_{\lb\in(0,1]}$ is equi-Lipschitz and equi-bounded. 
 \end{prop}
\proof
By (2) of Proposition \ref{prop:dom-fun} and item (5) of Lemma \ref{lem2:proT}, we derive $\{u^-_\lambda(x)\}_{\lb\in(0,1]}$ is Lipschitz continuous with the Lipschitz constant $\kappa$ independent of $\lambda$. Recall that $\{u^+_\lambda(x)\}_{\lb>0}$ is equi-bounded, namely
$$
-K_0\leq u^+_\lambda(x)\leq K_0,\quad \forall x\in M, \lambda>0
$$
for some constant $K_0$. Then (2) of Proposition \ref{prop:dom-fun} and Lemma \ref{lem2:Tmono} indicate
$$
u^-_\lambda(x)=\lim_{t\to+\infty}T^{\lambda,-}_t u^+_\lambda(x)\geq u^+_\lambda(x)\geq -K_0,\ \forall x\in M,\lambda\in(0,1].
$$
To show $\{u^-_\lambda(x)\}_{\lambda\in(0,1]}$ is bounded from above, it suffices to show $\{u^-_\lambda(x)|x\in\cG_\lb,\lambda\in(0,1]\}$ is bounded, since $\{u^-_\lambda\}_{\lambda\in(0,1]}$ is Lipschitzian with Lipschitz constant $\kappa$ independent of $\lambda\in(0,1]$. Indeed,
for each $x\in M,z\in \mathcal{G}_\lambda$,
$$
u^-_\lambda(x)\leq u^-_\lambda(z)+\kappa\cdot d(x,z).
$$
 Let $\gamma_{z,\lambda}:(-\infty,+\infty)\to M$ be a globally calibrated curve by $u^+_\lambda$ with $\gamma_{z,\lambda}(0)=z$. Then, for $t>0$,
 \begin{align*}
 	T^{\lambda,-}_tu^+_\lambda(z)&\leq u^+_\lambda(\gamma_{z,\lambda}(-t))e^{\lambda t}+\int^0_{-t}e^{-\lambda\tau}(L(\gamma_{z,\lambda}(\tau),\dot{\gamma}_{z,\lambda}(\tau))+c(H))\mbox{d}\tau\\
 	&=u^+_\lambda(z)\leq K_0.
 \end{align*}
We derive $u^-_\lambda(z)\leq K_0$ on $\mathcal{G}_\lambda$ for all $\lambda\in(0,1]$.
\qed

\section{Convergence of ground state solutions}\label{s3}

Before we prove Theorem \ref{thm:1}, we first propose the following definition:

\begin{defn}[Conjugate weak KAM solutions \cite{F}]\label{defn:conj-p}
A backward (resp. forward) $0-$weak KAM solution $u^-$ (resp. $u^+$) of \eqref{eq:hj} is said to be conjugated to $u^+$ (resp. $u^-$), if 
\begin{itemize}
\item $u^+\leq u^-$;
\item $u^-=u^+$ on the {\sf projected Mather set} $\cM$\footnote{see Definition \ref{defn:mat}}.
\end{itemize}
Such a pair $(u^-,u^+)$ is called {\sf conjugate}.
\end{defn}

\begin{lem}[Corollary 5.1.3. of \cite{F}]\label{lem:conj}
Any backward $0-$weak KAM solution $u^-$ of \eqref{eq:hj} conjugates to one and only one  forward $0-$weak KAM solution $u^+$, vice versa.
\end{lem}

On the other side, due to an analogue skill as in \cite{DFIZ}, we get:
\begin{lem}
The function $u_\lb^+$ in \eqref{eq:for-weak} converges to a uniquely identified forward $0-$weak KAM solution $u_0^+$ of \eqref{eq:hj}.
\end{lem}
\proof
From the proof of Proposition (\ref{prop:dom-fun}), $\widehat{u}_\lambda^-=-u^+_\lambda$ is the unique viscosity solution of (\ref{eq:sym-dis}). 
Due to Theorem 1.1 of \cite{DFIZ}, $\widehat{u}_\lambda^-$ converges to a unique backward $0$-weak KAM solution $\wh u_0^-$ of 
\[
\wh H(x,d_x u)=c(H),\quad x\in M.
\]
Due to Theorem 4.9.3 of \cite{F}, $u_0^+:=-\wh u_0^-$ is a forward $0-$weak KAM solution of 
\[
 H(x,d_x u)=c(H),\quad x\in M
\]
and $u_0^+=-\lim_{\lb\rightarrow 0_+}\wh u_\lb^-=\lim_{\lb\rightarrow 0_+} u_\lb^+$.
\qed

\vspace{20pt}

\noindent{\it Proof of Theorem \ref{thm:1}:} Suppose $u_0^-$ is the uniform limit of $u_{\lb_n}^-$ as $\lb_n\rightarrow 0_+$, then $u_0^-\geq u_0^+$. On the other side, $u_\lb^-=u_\lb^+$ on $\cG_\lb$, 
if \eqref{upper-aubry} holds, then for any Aubry class $\Lb_i$, there exists at least one point $z\in\Lb_i$ which can be accumulated by a sequence $\{z_n\in\cG_{\lb_n}\}_{n\in\N}$. Consequently, $u_0^-(z)=u_0^+(z)$. Due to Lemma \ref{lem:con-vis}, we get $u_0^-=u_0^+$ on $\Lb_i$ and further $u_0^-=u_0^+$ on $\cA$. Recall that $\cM\subset\cA$ due to Proposition \ref{prop:mat}, then $u_0^-$ is conjugated $u_0^+$, which is unique due to Lemma \ref{lem:conj}. So we get the uniqueness of $u_\lb^-$ as $\lb\rightarrow 0_+$.


For a generic $C^1$ Tonelli Hamiltonian $H(x,p)$, there exists a unique ergodic Mather measure (see Proposition \ref{prop:uni-erg}), then $\cA\slash d_c$ is a singleton due to Proposition \ref{prop:aub-cla} (namely $\cA$  consists of a unique Aubry class). Furthermore, by Proposition \ref{prop:u-semi} and Proposition \ref{prop:aub-min}, the condition \eqref{upper-aubry} always holds.
\qed

\vspace{40pt}

\appendix

\section{Aubry Mather Theory of Tonelli Hamiltonians}\label{a1}

As is known, the $C^2$ Tonelli Hamiltonian $H(x,p)$ has a dual {\sf Tonelli Lagrangian} 
\[
L(x,v):=\max_{p\in T_x^*M}\langle p,v\rangle-H(x,p),\quad (x,v)\in TM
\]
which is also $C^2$ and strictly convex, superlinear in $v$.  Consequently, for any $x,y\in M$ and $t>0$, the {\it action function}
\be\label{eq:act}
h^t(x,y):=\inf_{\substack{\gamma\in C^{ac}([0,t],M)\\\gamma(0)=x,\gamma(t)=y}}\int_0^tL(\gamma,\dot\gamma)+c(H)ds
\ee
always attains its infimum at a $C^2-$smooth minimizing curve $\gamma_{\min}:[0,t]\rightarrow M$, satisfying the {\sf Euler-Lagrange equation}
\be\label{eq:el}
\left\{
\begin{aligned}
&\frac{d}{dt}x=v,\\
&\frac d{dt}L_v(x,v)=L_x(x,v),
\end{aligned}
\right.
\ee
due to the {\sf Tonelli Theorem} and the {\sf Weierstrass Theorem}, see \cite{Mat}. A curve $\gamma:\R\rightarrow M$ is called {\sf critical}, if $(\gamma,\dot\gamma)$ solves (\ref{eq:el}). 
Denote the Lagrangian flow by $\Phi_{L}^t:TM\rightarrow TM$, then $\Phi_{L}^t$ is well defined for $t\in\R$ since $H(x,L_v(x,v))$ is invariant w.r.t. it.\medskip

\begin{defn} In \cite{Mat2}, the {\sf Peierls barrier} function 
\be
h^\infty(x,y):=\liminf_{t\rightarrow+\infty}h^t(x,y)
\ee
is proved to be well-defined and continuous on $M\times M$. Consequently, 
the {\sf projected Aubry set} is defined by 
\[
\cA:=\{x\in M: h^\infty(x,x)=0\}.
\]
With respect to the pseudo metric
\[
d_c(x,y):=h^\infty(x,y)+h^\infty(y,x),\quad\forall x,y\in\cA,
\]
we can decompose $\cA$ into a bunch of connected subsets ({\sf static classes} in \cite{CP}), such that any two point in the same class has a trivial $d_c-$distance. Without loss of generality, let's denote by $\cA/d_c$ the set of all the static classes.  

\end{defn}


Consider $TM$ (resp. $M$) as a measurable space and $\PP(TM)$ (resp. $\PP(M)$) by the set of all Borel probability measures on it.  A measure on $TM$ is denoted by $\wt\mu$, and we remove the tilde if we project it to $M$. We say that a sequence
$\{\wt\mu_n \}_n$ of probability measures weakly converges to a probability 
measure $\wt\mu$ if
\[
\lim_{n\rightarrow+\infty}\int_{TM}f(x,v)d\wt\mu_n(x,v)=\int_{TM} f(x,v)d\wt\mu(x,v)
\]
for any $f\in C_c(TM,\R)$. Accordingly, the deduced probability measure $\mu_n$ weakly converges to $\mu$, i.e.
\be\label{eq:pro-mes}
\lim_{n\rightarrow+\infty}\int_M f(x)d\mu_n(x)&:=&\lim_{n\rightarrow+\infty}\int_{TM} f\circ \pi(x,v) d\wt\mu_n(x,v)\nonumber\\
&=&\int_{TM} f\circ\pi(x,v)d\wt\mu(x,v)=:\int_M f(x)d\mu(x)
\ee
for any $f\in C(M,\R)$.

\begin{defn}
A probability measure $\wt\mu$ on $TM$ is {\sf closed} if it satisfies:
\begin{itemize}
\item $\int_{TM}|v|d\wt\mu(x,v)<+\infty$;
\item $\int_{TM}\langle \nabla\phi(x),v\rangle d\wt\mu(x,v)=0$ for every $\phi\in C^1(M,\R)$.
\end{itemize}
\end{defn}
Let's denote by $\PP_c(TM)$ the set of all closed measures on $TM$, then the following conclusion is proved in \cite{Mn1}:
\begin{thm}\label{thm:mane}
$\min_{\wt\mu\in\PP_c(TM)}\int_{TM}L(x,v)d\widetilde{\mu}(x,v)=-c(H)$. Moreover, the minimizer is $\Phi_{L}^t-$invariant and called a {\sf Mather measure}.
\end{thm}
\begin{defn}\label{defn:mat}
Define by $\wt{\mathfrak M}$ the set of Mather measures, which can be projected to $\mathfrak M\subset\PP(M)$ consisting of all the {\sf projected Mather measures} due to \eqref{eq:pro-mes}. The {\sf projected Mather set} is defined by 
\[
\cM:=\overline{\bigcup_{{\mu}\in{\mathfrak M}}supp({\mu})}\subset M.
\]
\end{defn} 
\begin{prop}\cite{B}\label{prop:mat}
$\cM\subset\cA$.
\end{prop}
\begin{defn}[Ma\~n\'e's genericity\cite{Mn1}]\label{defn:man-gene}
A property is called {\sf ($C^2-$)generic (in the sense of Ma\~n\'e)} for $H(x,p)$, if there exists a residual set $\mathcal O\subset C^2(M,\R)$ such that for any $\psi\in \mathcal O$, the property holds for $H+\psi$. Accordingly, a Tonelli Hamiltonian $H(x,p)$ is called {\sf generic}, if we can find another Tonelli $H_0(x,p)$ and a residue set $\cO\subset C^2(M,\R)$, such that $H-H_0\in \cO$.
\end{defn}
\begin{prop}[Theorem C of \cite{Mn1}]\label{prop:uni-erg}
For a generic $C^2$ Tonelli Hamiltonian, the associated Mather measure is uniquely ergodic. 
\end{prop}
\begin{prop}[Lemma 5.3 of \cite{CP}]\label{prop:aub-cla}
If the Mather measure is uniquely ergodic, the Aubry class has to be unique. 
\end{prop}
\begin{prop}[Proposition 5.3 of \cite{B}]\label{prop:aub-min}
For autonomous Tonelli Lagrangian, if $\cA$ is of a unique Aubry class, then $\cA=\cG$ with 
\ben
\mathcal{G}=\big\{x\in M\big| \text{there exists }\gamma:\mathbb{R}\to M\  \mbox{with}\  \gamma(0)=x, \text{such that }\\
\forall a<b\in\R, \gamma|_{t\in[a,b]} \text{ realizes } h^{b-a}(\gamma(a),\gamma(b))\big\}.
\een
\end{prop}

\section{Viscosity solutions of \eqref{eq:hj}}\label{a2}

\begin{lem}[Theorem 7.6.2. of \cite{F}]\label{lem:sta-weak}
The backward $0-$weak KAM solution has to be a viscosity solution, vice versa.
\end{lem}
\begin{lem}[item 3 of Remark 4.9.3 in \cite{C}]
\begin{itemize}
\item For any $y\in M$ fixed, $h^\infty(y,\cdot)$ is a backward $0-$weak KAM solution.
\item for any $y\in M$ fixed, $-h^\infty(\cdot,y)$ is a forward $0-$weak KAM solution;
\end{itemize}
\end{lem}

\begin{prop}[Theorem 8.6.1 of \cite{F}]\label{prop:exp-weak}
Any viscosity solution $u$ of \eqref{eq:hj} can be formally expressed by 
\[
u(x):=\inf_{x_0\in\cA}\{u(x_0)+h^\infty(x_0,x)\},\quad\forall x\in M.
\]
\end{prop}
\begin{lem}\cite{FS}\label{lem:con-vis}
For any $\Lb_i\in\cA/d_c$, any two viscosity solutions of \eqref{eq:hj} differs by a constant on $\Lambda_i$.
\end{lem}
\proof
This conclusion is a direct corollary of previous Proposition \ref{prop:exp-weak}. Precisely, for any viscosity solution $u$, we have 
\ben
u(x)&=&\inf_{x_0\in\cA}\{u(x_0)+h^\infty(x_0,x)\}\\
&=&\inf_{x_0\in\cA}\{u(x_0)+h^\infty(x_0,y)+h^\infty(y,x)\}\\
&=&\inf_{x_0\in\cA}\{u(x_0)+h^\infty(x_0,y)\}+h^\infty(y,x)\\
&=& u(y)+h^\infty(y,x)
\een
as long as $x,y$ belonging to the same static class. Therefore, 
\[
\om(y)-u(y)=\om(x)-u(x),\quad\forall x,y\in\Lb_i
\]
for any two viscosity solutions $u$ and $\om$.
\qed


\section{Properties of $T_{t}^{\lb,-}$}

\begin{lem}\label{lem2:proT}
\begin{description}
	\item [(1).] For $s,t>0$ and $\phi(x)\in C(M,\mathbb{R})$,
	$$
	T^{\lambda,-}_{t+s}\phi(x)=T^{\lambda,-}_t\circ T^{\lambda,-}_s\phi(x).
	$$
	\item[(2).] Let $\phi_1(x),\phi_2(x)\in C(M,\mathbb{R})$. Then,
	$$
	\|T^{\lambda,-}_t\phi_1-T^{\lambda,-}_t\phi_2\|\leq e^{\lambda t}\|\phi_1-\phi_2\|,\quad\forall t>0.
	$$
\item[(3).] For $\lambda\in(0,1]$, each minimizer $\gamma$ of $T^{\lambda,-}_t\phi(x)$ is $C^2$ and 
$|\dot{\gamma}(\tau)|\leq \alpha_0$, where $\alpha_0$ is independent of $\lambda$.
\item [(4).] Let $\phi_1,\phi_2\in C(M,\mathbb{R})$ and $\phi_1(x)\leq\phi_2(x),x\in M$. Then,
\begin{equation}
	T^{\lambda,-}_t\phi_1(x)\leq T^{\lambda,-}_t\phi_2(x), \quad \forall t\geq0, x\in M.
\end{equation}	
\item[(5).] For $\lambda\in(0,1]$ and $t>\mbox{diam}(M)$, the map $x\longmapsto T^{\lambda,-}_t\phi(x)$ is equi-Lipschitz, i.e.,
	$$
	|T^{\lambda,-}_t\phi(x)-T^{\lambda,-}_t\phi(y)|\leq \kappa d(x,y), \quad t>\mbox{diam}(M),
	$$
	where $\kappa$ is independent of $\lambda,\phi$ and $t$.
\item[(6).] $u\in C(M,\mathbb{R})$ is a backward $\lambda$-weak KAM solution if and only if 
$$
T^{\lambda,-}_tu(x)=u(x),\forall x\in M,t>0.
$$
\item[(7).] Each fixed point of $T_t^{\lb,-}$ is a viscosity solution of \eqref{eq:dis}, vice versa.
\end{description}	
\end{lem}
\proof

 (1). Note that 
 $$
 h^{-(s+t),0}(y,x)=\inf_{z\in M}\{h^{-(s+t),-t}(y,z)+h^{-t,0}(z,x)\}
 $$ 
 and 
 $$
 h^{-(s+t),-t}(y,z)=e^{\lambda t}h^{-s,0}(y,z).
 $$ 
 We derive
\begin{align*}
	T^{\lambda,-}_{s+t}\phi(x)&=\inf_{y\in M}\{\phi(y)e^{\lambda(s+t)}+h^{-(s+t),0}(y,x)\}\\
	&=\inf_{y\in M}\{\phi(y)e^{\lambda(s+t)}+\inf_{z\in M}\{h^{-(s+t),-t}(y,z)+h^{-t,0}(z,x)\}\}\\
	&=\inf_{z\in M}\{\inf_{y\in M}\{\phi(y)e^{\lambda(s+t)}+e^{\lambda t}h^{-s,0}(y,z)\}+h^{-t,0}(z,x)\}\\
	&=\inf_{z\in M}\{e^{\lambda t}T^{\lambda,-}_{s}\phi(z)+h^{-t,0}(z,x)\}\\
	&=T^{\lambda,-}_t\circ T^{\lambda,-}_s\phi(x).
\end{align*}

(2). For every $x\in M$, let $\gamma_1:[-t,0]\to M$ be a minimizer of $T^{\lambda,-}_{t}\phi_1(x)$. Then,
\begin{align*}
	T^{\lambda,-}_t\phi_2(x)-T^{\lambda,-}_t\phi_1(x)&\leq \phi(\gamma_1(-t))e^{\lambda t}-\phi(\gamma_2(-t))e^{\lambda t}\\
	&\leq e^{\lambda t}\|\phi_2-\phi_1\|.
\end{align*}
Similarly, we can obtain
$$
T^{\lambda,-}_t\phi_1(x)-T^{\lambda,-}_t\phi_2(x)\leq e^{\lambda t}\|\phi_2-\phi_1\|.
$$
Hence, $$
\|T^{\lambda,-}_t\phi_1-T^{\lambda,-}_t\phi_2\|\leq e^{\lambda t}\|\phi_1-\phi_2\|.
$$

(3). By the Weierstrass Theorem, we derive that $\gamma$ solves Euler-Lagrangian equation (\ref{eq2:E-L}). Hence, $\gamma$ is $C^2$. By the definition of $T^{\lambda,-}_t\phi(x)$, for $\frac{1}{2}<s_2-s_1<1$
$$
e^{\lambda s_1}T^{\lambda,-}_{-s_1+t}\phi(\gamma(-s_1))-e^{\lambda s_2}T^{\lambda,-}_{-s_2+t}\phi(\gamma(-s_2))=\int^{-s_1}_{-s_2}e^{-\lambda\tau}(L(\gamma(\tau),\dot{\gamma}(\tau))+c(H))\mbox{d}\tau.
$$
On the other hand, let $\beta:[-s_2,-s_1]\to M$ be a geodesic satisfying
$\beta(-s_1)=\gamma(-s_1),\beta(-s_2)=\gamma(-s_2)$ and $|\dot{\gamma}(\tau)|\leq 2\mbox{diam}(M)=:k_0$.
\begin{align*}
	e^{\lambda s_1}T^{\lambda,-}_{-s_1+t}\phi(\gamma(-s_1))-e^{\lambda s_2}T^{\lambda,-}_{-s_2+t}\phi(\gamma(-s_2))&\leq\int^{-s_1}_{-s_2}e^{-\lambda\tau}(L(\beta(\tau),\dot{\beta}(\tau))+c(H))\mbox{d}\tau\\
	&\leq (C_{k_0}+c(H))\int^{-s_1}_{-s_2}e^{-\lambda \tau}\mbox{d}\tau.
\end{align*}
Then,
\begin{align*}
	\int^{-s_1}_{-s_2}e^{-\lambda \tau}(|\dot{\gamma}(\tau)|-C(1)+c(H))\mbox{d}&\tau\leq \int^{-s_1}_{-s_2}e^{-\lambda\tau}(L(\gamma(\tau),\dot{\gamma}(\tau))+c(H))\mbox{d}\tau\\
	&\leq (C_{k_0}+c(H))\int^{-s_1}_{-s_2}e^{-\lambda \tau}\mbox{d}\tau.
\end{align*}
There exists $\tau_0\in(-s_2,-s_1)$ such that
$|\dot{\gamma}(\tau_0)|\leq C_{k_0}+C(1)$.
Note that $\gamma(\tau)$ solves Euler-Lagrangian equation
(\ref{eq2:E-L}). We derive $|\dot{\gamma}(\tau)|\leq \alpha_0$.\medskip

(4). For each $\gamma\in C^{ac}([-t,0],M)$ with $\gamma(0)=x$, it holds
\begin{align*}
	&\phi_1(\gamma(-t))e^{\lambda t}+\int^0_{-t}e^{-\lambda\tau}(L(\gamma(\tau),\dot{\gamma}(\tau))+c(H))\mbox{d}\tau\\
	&\leq \phi_2(\gamma(-t))e^{\lambda t}+\int^0_{-t}e^{-\lambda\tau}(L(\gamma(\tau),\dot{\gamma}(\tau))+c(H))\mbox{d}\tau.
\end{align*}
Then, $T^{\lambda,-}_t\phi_1(x)\leq T^{\lambda,-}_t\phi_2(x)$.\medskip

(5). Let $\gamma_x:[-t,0]\to M$ be a minimizer of $T^{\lambda,-}_t\phi(x)$ and $\Delta t=d(x,y)$ and let $\beta:[-\Delta t,0]\to M$ be a geodesic satisfying $\beta(-\Delta t)=\gamma_x(-\Delta t), \beta(0)=y$, and
$$
|\dot{\beta}(\tau)|\equiv \frac{d(\gamma_x(-\Delta t),y)}{\Delta t}\leq 
\frac{d(\gamma_x(-\Delta t),x)}{\Delta t}+1\leq \alpha_0+1,
$$
where $|\dot{\gamma}_x(\tau)|\leq \alpha_0$.
\begin{align*}
	T^{\lambda,-}_t\phi(y)-T^{\lambda,-}_t\phi(x)&\leq \int^0_{-\Delta t}e^{-\lambda \tau}\big(L(\beta(\tau),\dot{\beta}(\tau))-L(\gamma(\tau),\dot{\gamma}(\tau))\big)\mbox{d}\tau\\
	&\leq (C_{\alpha_0+1}+C(0))\int^0_{-\Delta t}e^{-\lambda\tau}\mbox{d}\tau\\
	&=(C_{\alpha_0+1}+C(0))\cdot d(x,y).
\end{align*}
Similarly, we have
$$
T^{\lambda,-}_t\phi(x)-T^{\lambda,-}_t\phi(y)\leq (C_{\alpha_0+1}+C(0))\cdot d(x,y).
$$
Let $\kappa=C_{\alpha_0+1}+C(0)$. We have
$$
|T^{\lambda,-}_t\phi(y)-T^{\lambda,-}_t\phi(x)|\leq \kappa\cdot d(x,y).
$$

(6). Let $u(x)$ be a backward $\lambda$-weak KAM solution and $\gamma_{x,\lambda}^-:(-\infty,0]\to M$ be a backward calibrated curve satisfying $\gamma_{x,\lambda}^-(0)=x$.
For each $t>0$,
$$
u(x)-e^{\lambda t}u(\gamma_{x,\lambda}^-(-t))=\int^0_{-t}e^{-\lambda\tau}(L(\gamma_{x,\lambda}^-(\tau),\dot{\gamma}_{x,\lambda}^-(\tau))+c(H))\mbox{d}\tau.
$$
By $u\prec_\lambda L+c(H)$, we derive 
$$
u(x)=\inf_{\substack{ \gamma\in C^{ac}([-t,0],M)\\\gamma(0)=x}}\bigg\{e^{\lambda t}u(\gamma(-t))+\int^0_{-t}e^{-\lambda\tau}(L(\gamma(\tau),\dot{\gamma}(\tau))+c(H))\mbox{d}\tau\bigg\}=T^{\lambda,-}_tu(x).
$$

On the other hand, we assume $T^{\lambda,-}_tu(x)=u(x)$.
Let $\gamma\in C^{ac}([t_1,t_2],M)$. Then,
$$
T^{\lambda,-}_{t_2-t_1}u(\tilde{\gamma}(0))\leq e^{\lambda(t_2-t_1)}u(\tilde{\gamma}(t_1-t_2))+\int^0_{t_1-t_2}e^{-\lambda \tau}(L(\wt{\gamma}(\tau),\dot{\wt{\gamma}}(\tau))+c(H))\mbox{d}\tau,
$$
where $\tilde{\gamma}\in C^{ac}([t_1-t_2,0],M)$ is defined by
$\tilde{\gamma}(s)=\gamma(s+t_2)$.
From $T^{\lambda,-}_{t_2-t_1}u(x)=u(x)$, it follows that
$$
u(\gamma(t_2))\leq e^{\lambda (t_2-t_1)}u(\gamma(t_1))+e^{\lambda t_2}\int^{t_2}_{t_1}e^{-\lambda\tau}(L(\gamma(\tau),\dot{\gamma}(\tau))+c(H))\mbox{d}\tau,
$$
which implies
$$
e^{-\lambda t_2}u(\gamma(t_2))\leq e^{-\lambda t_1}u(\gamma(t_1))+\int^{t_2}_{t_1}e^{-\lambda\tau}(L(\gamma(\tau),\dot{\gamma}(\tau))+c(H))\mbox{d}\tau.
$$
Hence, $u\prec_\lambda L+c(H)$.

For each $n\in\mathbb{N}$, we assume $\gamma_n:[-n,0]\to M$ is a minimizer of $T^{\lambda,-}_nu(x)$. Then, for each $t\in[0,n]$, 
\begin{align*}
u(x)-e^{\lambda t}u(\gamma_n(-t))&=T^{\lambda,-}_nu(x)-e^{\lambda t}T^{\lambda,-}_{n-t}u(\gamma_n(-t))\\
&=\int^0_{-t}e^{-\lambda\tau}(L(\gamma_n(\tau),\dot{\gamma}_n(\tau))+c(H))\mbox{d}\tau.
\end{align*}
Note that $|\dot{\gamma}_n(\tau)|\leq \alpha_0$ for all $n\in M$. By the Ascoli Theorem, there exists a subsequence $\{\gamma_{n_k}\}$, such that
$\gamma_{n_k}$ converges uniformly to $\gamma_x\in C^{ac}((-\infty,0],M)$ on any finite interval of $(-\infty,0]$.
Then, for each $t>0$, we derive
$$
u(x)-e^{\lambda t}u(\gamma_x(-t))\geq\int^0_{-t}e^{-\lambda\tau}(L(\gamma_x(\tau),\dot{\gamma}_x(\tau))+c(H))\mbox{d}\tau,
$$
which implies
$$
u(x)-e^{\lambda t}u(\gamma_x(-t))=\int^0_{-t}e^{-\lambda\tau}(L(\gamma_x(\tau),\dot{\gamma}_x(\tau))+c(H))\mbox{d}\tau,
$$
since $u\prec_\lambda L+c(H)$.
This means $\gamma_x:(-\infty,0]\to M$ is a backward calibrated curve. Hence, $u$ is a $\lambda$-weak KAM solution.\medskip

(7). Assume $u(x)=T^{\lambda,-}_tu(x)$ for $t>0$. Due to (6), $u(x)$ is a backward $\lambda$ -weak KAM solution of (\ref{eq:dis}).
Let $x_0\in M, v_0\in T_{x_0}M$ and let $\phi(x)\in C^1(M,\mathbb{R})$ and $u-\phi$ attains maximum at $x_0$. For $\Delta t<0$, we assume $\gamma:[t+\Delta t,t]\to M$ is an absolutely continuous curve with $\gamma(t)=x_0$ and $\dot{\gamma}(t)=v_0$. Then,
\begin{align*}
	e^{-\lambda t}\big(\phi(\gamma(t))-\phi(\gamma(t+\Delta t))\big)&\leq 
	e^{-\lambda t}u(\gamma(t))-e^{-\lambda t}u(\gamma(t+\Delta t))\\
	&\leq e^{-\lambda t}u(\gamma(t))-e^{-\lambda(t+\Delta t)}u(\gamma(t+\Delta t))\\&+e^{-\lambda(t+\Delta t)}u(\gamma(t+\Delta t))-e^{-\lambda t}u(\gamma(t+\Delta t)).
\end{align*}
Then,
\begin{align*}
	e^{-\lambda t}\frac{\phi(\gamma(t+\Delta t))-\phi(\gamma(t))}{\Delta t}&\leq 
	\frac{1}{\Delta t}\int^{t+\Delta t}_te^{-\lambda\tau}(L(\gamma(\tau),\dot{\gamma}(\tau))+c(H))\mbox{d}\tau\\
	&-\frac{e^{-\lambda(t+\Delta t)}-e^{-\lambda t}}{\Delta t}\cdot u(\gamma(t+\Delta t)).
\end{align*}
Taking $\Delta t\to 0^-$, we derive 
$$
e^{-\lambda t}d_x\phi(x_0)\cdot v_0\leq e^{-\lambda t}(L(x_0,v_0)+c(H))+\lambda e^{-\lambda t}\cdot u(x_0),
$$
which implies
$$
-\lambda u(x_0)+d_x\phi(x_0)\cdot v_0-L(x_0,v_0)\leq c(H).
$$
Hence,
$$
-\lambda u(x_0)+H(x_0,d_x\phi(x_0))\leq c(H).
$$
On the other hand, let $\psi(x)\in C^1(M,\mathbb{R})$ and $u-\psi$ attains the minimum at $x_0$ and let $\gamma_x:[t+\Delta t,t]\to M$ be a calibrated curve by $u$ with $\gamma_x(t)=x_0$. Then,
$$
\psi(\gamma_x(t))-\psi(\gamma_x(t+\Delta t))\geq u(\gamma_x(t))-u(\gamma_x(t+\Delta t)).
$$
Note that
\begin{align*}
	&e^{-\lambda t}\big(u(\gamma_x(t))-u(\gamma_x(t+\Delta t))\big)\\
&=e^{-\lambda t}u(\gamma_x(t))-e^{-\lambda(t+\Delta t)}u(\gamma_x(t+\Delta t))\\
&\hspace{1cm}+e^{-\lambda(t+\Delta t)}u(\gamma_x(t+\Delta t))-e^{-\lambda t}u(\gamma_x(t+\Delta t))\\
&=\int^t_{t+\Delta t}e^{-\lambda\tau}(L(\gamma_x(\tau),\dot{\gamma}_x(\tau))+c(H))\mbox{d}\tau+(e^{-\lambda(t+\Delta t)}-e^{-\lambda t})u(\gamma_x(t+\Delta t)).
\end{align*}
We derive 
\begin{align*}
	e^{-\lambda t}\cdot\frac{\psi(\gamma_x(t+\Delta t))-\psi(\gamma_x(t))}{\Delta t}&\geq \frac{1}{\Delta t}\int^{t+\Delta t}_te^{-\lambda\tau}(L(\gamma_x(\tau),\dot{\gamma}_x(\tau))+c(H))\mbox{d}\tau\\
	&-\frac{e^{-\lambda(t+\Delta t)}-e^{-\lambda t}}{\Delta t}\cdot u(\gamma_x(t+\Delta t)).
\end{align*}
Taking $\Delta t\to 0^-$, we derive that
$$
-\lambda u(x_0)+d_x\psi(x_0)\cdot v_0-L(x_0,v_0)\geq c(H),
$$
which implies
$$
-\lambda u(x_0)+H(x_0,d_x\psi(x_0))\geq c(H).
$$
so $u$ is a viscosity solution of \eqref{eq:dis}.

Suppose $\omega(x)$ is a viscosity solution of \eqref{eq:dis}, then 
$\omega(x)$ is 
Lipschitz due to the superlinearity of $H(x,p)$, see \cite{Ba}. For the reduced Lipschitz Lagrangian $\mathbf L^\lb(x,v):=L(x,v)+\lb \om(x)$, $\om(x)$ is also the viscosity solution of $\mathbf{H}^\lambda(x,d_x\omega(x))=c(H)$, where 
\[
\mathbf{H}^\lambda(x,p)=H(x,p)-\lambda \omega(x)
\]
 is the corresponding Hamiltonian. Then, $U:(x,t)\in M\times[0,+\infty)\to M$ defined by
$$
U(x,t):=\inf_{\substack{\gamma\in C^{ac}([0,t],M) \\\gamma(0)=x }}\{\omega(\gamma(-t))+\int^0_{-t}\mathbf{L}^\lambda(\gamma(\tau),\dot{\gamma}(\tau))+c(H)\mbox{d}\tau\}, \quad\forall t\geq 0, 
$$
is a viscosity solution of the Cauchy problem 
 \begin{equation}\label{eq_app}
 	\left\{
 \begin{aligned}
 &\partial_tu+\mathbf{H}^\lambda(x,d_x\omega)=c(H),\\
 &u(x,0)=\om(x),\quad t\geq 0.
 \end{aligned}
 \right.
 \end{equation}

  Note that $\omega(x)$ is also a viscosity of $\mathbf{H}^\lambda(x,d_x\omega(x))=c(H)$. We derive $\omega(x)$ is a solution to the Cauchy problem (\ref{eq_app}). From the uniqueness of viscosity solution, it follows that $U(x,t)=\omega(x)$ for $x\in M,t\geq 0$. Hence, for each absolutely continuous curve $\gamma:[s,t]\to M$,
  $$
  \omega(\gamma(t))-\omega(\gamma(s))\leq\int^t_s\mathbf{L}^\lambda(\gamma,\dot\gamma)+c(H)\mbox{d}\tau.
  $$

Fix a sequence $\{t_n\}_{n\in\mathbb{N}}$ tending to $+\infty$ as $n\to\infty$.  Due to $U(x,t_n)=\omega(x)$, for each $n\in\mathbb{N}$, there exists an absolutely continuous curve $\gamma_n:[-t_n,0]\to M$ such that $\gamma_n(0)=x$ and 
$$
\omega(x)=\omega(\gamma_n(-t_n))+\int^0_{-t_n}\mathbf{L}^{\lambda}(\gamma_n(\tau),\dot{\gamma}_n(\tau))+c(H)\mbox{d}\tau.
$$ 


From the superlinearity of $\mathbf{L}^\lambda$ in $v$ and Lipschitz continuity of $\omega$, we derive $\{\|\dot{\gamma}_n\|_{L^\infty}\}_n$ is equi-bounded. By Ascoli Theorem, there exists an 
 subsequence of $\{\gamma_n\}$ (denoted still by $\gamma_n$) uniformly converging to an
 absolutely continuous curve $\gamma_*:(-\infty,0]\to M$ on each finite interval of $(-\infty,0]$, such that $\gamma_*(0)=x$ and 
 \begin{equation}\label{eqD:calibrated}
 	\omega(x)-\omega(\gamma_*(-t))=\int^0_{-t}\mathbf{L}^\lambda(\gamma_*(\tau),\dot{\gamma}_*(\tau))+c(H)\mbox{d}\tau, t>0. 
 \end{equation}
 Let $\gamma\in C^{ac}([a,b], M)$. Then, $\omega(\gamma(\tau))$ and 
 $$
 s\longmapsto \int^0_sL(\gamma(\tau),\dot{\gamma}(\tau))+c(H)+\lambda\omega(\gamma(\tau))\mbox{d}\tau
 $$
 are differentiable a.e. on $[a,b]$.
 For $t\in [a,b]$ and $\Delta t\not=0$ with $t+\Delta t\in[a,b]$,
 $$
 \frac{\omega(\gamma(t+\Delta t))-\omega(\gamma(t))}{\Delta t}\leq\frac{1}{\Delta t}\int^{t+\Delta t}_t L(\gamma(\tau),\dot{\gamma}(\tau))+c(H)+\lambda\omega(\gamma(\tau))\mbox{d}\tau.
 $$ 
  Taking $\Delta t$ tending to $0$, we derive that
 $$
 \frac{\mbox{d}\omega(\gamma(t))}{\mbox{d}t}\leq L(\gamma(t),\dot{\gamma}(t))+c(H)+\lambda\omega(\gamma(t)), a.e.\  t\in [a,b].
 $$
 Then,
 $$
 \frac{\mbox{d}}{\mbox{d}t}\bigg(e^{-\lambda t}\omega(\gamma(t))\bigg)\leq e^{-\lambda t}\big(L(\gamma(t),\dot{\gamma}(t))+c(H)\big), a.e.\  t\in [a,b].
 $$
 Integrating on $[a,b]$, we derive
 $$
 e^{-\lambda b}\omega(\gamma(b))-e^{-\lambda a}\omega(\gamma(a))\leq \int^b_ae^{-\lambda\tau}\big(L(\gamma(\tau),\dot{\gamma}(\tau))+c(H)\big)\mbox{d}\tau,
 $$
 which implies $\omega\prec_\lambda L+c(H)$.
 By a similar discussion, we derive from (\ref{eqD:calibrated}) that $\gamma_*$ is a calibrated curve by $\omega$, i.e.,
 $$
 \omega(x)-e^{\lambda t}\omega(\gamma_*(-t))=\int^0_{-t}e^{-\lambda\tau}\big(L(\gamma_*(\tau),\dot{\gamma}_*(\tau))+c(H)\big)\mbox{d}\tau,\quad \forall t>0,
 $$
which implies $\om$ is a backward $\lb-$weak KAM solution of \eqref{eq:dis}.\qed

\section{Upper semi-continuity of $\cG_\lb$}\label{a3}

%

\begin{prop}[Upper semicontinuity]\label{prop:u-semi}
As a set-valued function,  
\[
\varlimsup_{\lb\rightarrow 0_+} \cG_\lb\subset \cG\subset M.
\] 
\end{prop}
\proof
Let $\gamma_n$ be a globally calibrated curve by $u^+_{\lambda_n}$ with the parameter $\lambda_n\to 0_+$, as $n\to\infty$. 
 To show the proposition, it suffices to show any accumulating curve $\gamma^*$ of $\{\gamma_n\}$  realizes $h^{b-a}(\gamma^*(a),\gamma^*(b)),a<b\in\mathbb{R}$.
 
 Otherwise, there exists an interval $[a,b]$ and a curve $\eta^*\in C^{ac}([a,b],M)$ such that $\eta^*(a)=\gamma^*(a),\eta^*(b)=\gamma^*(b)$, and 
 \begin{equation}\label{eqD:min}
 	h^{b-a}(\gamma^*(a),\gamma^*(b))=\int^b_aL(\eta^*,\dot{\eta}^*)+c(H)\mbox{d}\tau
  <\int^b_aL(\gamma^*,\dot{\gamma}^*)+c(H)\mbox{d}\tau.
 \end{equation}
 By Weierstrass Theorem,
 one can easily check $\eta^*$ is $C^2$ and $|\dot{\eta^*}|\leq \kappa_0$.
 For sufficiently large $n\in\mathbb{N}$, we define $\eta_n\in C^{ac}([a,b],M)$ by
 $$
 \eta_n(s)=
 \begin{cases}
 \beta_{1,n}(s), s\in [a,a+d_{1,n}], \\
 	\eta^*(s),s\in[a+d_{1,n},b-d_{2,n}], \\
 	\beta_{2,n}(s),s\in[b-d_{2,n},b],
 \end{cases}
 $$
 where $d_{1,n}=d(\gamma_n(a),\gamma^*(a)),d_{2,n}=d(\gamma_n(b),\gamma^*(b))$, $\beta_{1,n}$ is the geodesic connecting $\gamma_n(a)$ and $\eta^*(a+d_{1,n})$ with $|\dot{\beta}_{n}|\leq \kappa_0+1$, $\beta_{2,n}$ is the geodesic connecting $\eta^*(b-d_{2,n})$ and $\gamma_n(b)$ with $|\dot{\beta}_{2,n}|\leq \kappa_0+1$.
Then,
$\eta_n$ converges uniformly to $\eta^*$ on $[a,b]$ and 
$
\int^b_a|e^{-\lambda_n\tau}(L(\eta_n,\dot{\eta}_n)+c(H))|\mbox{d}\tau$ is bounded.

By Dominated Convergence Theorem, we derive
\begin{equation*}
	\lim_{n\to+\infty}\int^b_ae^{-\lambda_n\tau}(L(\eta_n,\dot{\eta}_n)+c(H))\mbox{d}\tau=\int^b_aL(\eta^*,\dot{\eta}^*)+c(H)\mbox{d}\tau
\end{equation*}
Hence,
\begin{align*}
	h^{b-a}(\gamma^*(a),\gamma^*(b))=\int^b_aL(\eta^*,\dot{\eta}^*)+c(H)\mbox{d}\tau &= \lim_{n\to+\infty}\int^b_ae^{-\lambda_n\tau}(L(\eta_n,\dot{\eta}_n)+c(H))\mbox{d}\tau\\
	&\geq \liminf_{n\to+\infty}\int^b_ae^{-\lambda_n\tau}(L(\gamma_n,\dot{\gamma}_n)+c(H))\mbox{d}\tau\\
	&\geq \int^b_aL(\gamma^*,\dot{\gamma}^*)+c(H)\mbox{d}\tau.\\
\end{align*}
Combining (\ref{eqD:min}), we derive a contradiction. Hence, for each $a<b\in\mathbb{R}$,
$$
h^{b-a}(\gamma^*(a),\gamma^*(b))=\int^b_aL(\gamma^*,\dot{\gamma}^*)+c(H)\mbox{d}\tau
$$
Then the assertion follows.
%
\qed


\vspace{40pt}

\begin{thebibliography}{}
\bibitem{Ba}  G. Barles. {\it Solutions de viscosit\'e des \'equations de Hamilton-Jacobi}, volume. 17 of {\it Math\'ematiques \& Applications}. Springer-Verlag, Paris, 1994.
\bibitem{B} P. Bernard. Connecting orbits of time dependent Lagrangian systems. {\it Ann. Inst. Fourier, Grenoble}, 52(5): 1533-1568, 2002.
\bibitem{CCIZ} Q. Chen, W. Cheng, H. Ishii, and K. Zhao. Vanishing contact structure problem and convergence of the viscosity solutions. {\it Comm. Partial Differential Equations}, 44(9): 801-836, 2019.

\bibitem{CIPP}G. Contreras, R. Iturriaga, G. P. Paternain, and M. Paternain. Lagrangian graphs, minimizing measures and Mane critical values. {\it Geom. Funct. Anal.}, 8: 788-809, 1998.

\bibitem{C} G. Contreras. {\it Global minimizers of autonomous Lagrangians}. 22 Col\'oquio Brasileiro de Matem\'atica, IMPA, Rio de Janeiro, 1999.
\bibitem{CP} G. Contreras and Gabriel P. Paternain. Connecting orbits between static classes for generic Lagrangian systems. {\it Topology}, 41(4): 645-666, 2002.
\bibitem{DFIZ} A. Davini, A. Fathi, R. Iturriaga, and M. Zavidovique. Convergence of the solutions of the discounted Hamilton-Jacobi equation: convergence of the discounted solutions. {\it Invent. Math.}, 206(1): 29-55, 2016.





\bibitem{DW} A. Davini and L. Wang. On the vanishing discount problem from the negative direction. {\it Discrete Contin. Dyn. Syst.}, 41(5): 2377-2389, 2021.
\bibitem{F} A. Fathi. {\it Weak KAM theorem in Lagrangian dynamics}. preliminary version 10, Lyon. unpublishied. 2008.
\bibitem{FS} A. Fathi and A. Siconolfi. PDE aspects of Aubry-Mather theory for quasiconvex Hamiltonians. {\it Calc. Var. Partial Differen. Equ}. 22: 185-228, 2005.
\bibitem{LPV} P. L. Lions, G. Papanicolaou, and S. Varadhan. {\it Homogenization of Hamilton-Jacobi equation}. unpublished preprint, 1987.
\bibitem{Mn1} R. Ma\~n\'e. Generic properties and problems of minimizing measures of Lagrangian systems. {\it Nonlinearity}, 9(2): 273-310, 1996. 
\bibitem{MS}  S. Mar\`o and A. Sorrentino. Aubry-Mather theory for conformally symplectic systems. {\it Comm. Math. Phys.}, 354: 775-808, 2017.
\bibitem{Mat} J. Mather. Action minimizing invariant measures for positive definite Lagrangian systems. {\it Math. Z.}, 207: 169-207, 1991.
\bibitem{Mat2} J. Mather. Variational construction of connecting orbits. {\it Ann. Inst. Fourier (Grenoble)}, 43(5): 1349-1386, 1993.
\bibitem{WWYNON} K. Wang, L. Wang, and J. Yan. Implicit variational principle for contact Hamiltonian systems. {\it Nonlinearity}, 30: 492-515, 2017.

\bibitem{WWYJMPA} K. Wang, L. Wang, and J. Yan. Variational principle for contact Hamiltonian systems and its applications. {\it J. Math. Pures Appl}., 123: 167-200, 2019.

\bibitem{WWY}  K. Wang, L. Wang, and J. Yan. Aubry-Mather theory for contact Hamiltonian systems, {\it Comm. Math. Phys}., 366: 981-1023, 2019.

\bibitem{WWY2} K. Wang, L. Wang, and J. Yan. Weak KAM solutions of Hamilton-Jacobi equations with decreasing dependence on unknown functions, {\it J. Differential Equations}, 286: 411-432, 2021.

\bibitem{WYZ} Y.-N. Wang, J. Yan, and J. Zhang. {\it Convergence of viscosity solutions of generalized contact Hamilton-Jacobi equations}. Arxiv:2004.12269, 2020.
\bibitem{Z} M. Zavidovique. {\it Convergence of solutions for some degenerate discounted Hamilton-Jacobi equations}, arXiv:2006.00779, 2020.
\end{thebibliography}
\end{document}